\documentclass[a4paper,11pt,french]{article}
\usepackage{amsfonts}
\usepackage[T1]{fontenc}
\usepackage[ansinew]{inputenc}
\usepackage{babel}
\usepackage[top=2.5cm,bottom=2.5cm,right=2cm,left=2cm]{geometry}

\newcommand{\N}{\rm I\!N}
\newcommand{\Z}{\mathbb {Z}}

\newcommand{\R}{\rm I\!R}
\newcommand{\C}{\mathbb{C}}
\newcommand{\T}{\mathbb{T}}

\begin{document}

\centerline{\bf \textsc{\Large{\'Equation cohomologique d'un }}}
\medskip
\centerline{\bf  \textsc{\Large{automorphisme affine hyperbolique
du tore}}}

\vskip 0.5cm

\centerline{Abdellatif ZEGGAR}

\vskip 0.2cm

\centerline{(Février 2019)}

\vskip 0.5cm

\noindent {\bf Résumé} :  {\it Soit $A$ une matrice hyperbolique
(matrice sans valeur propre de module égal à $1$) appartenant à
$\hbox{GL}(p,\Z)$  et soit $b$ un élément de $\R^p$.
L'automorphisme affine $\gamma : x\mapsto \gamma (x)=Ax+b$ de
$\R^p$ induit sur le tore $\displaystyle \mathbb{T}^p=\R^p/\Z^p$
un automorphisme affine hyperbolique (d'Anosov) que l'on note
encore $\gamma $.  L'objet de ce travail est l'étude de l'équation
cohomologique :
$$\displaystyle f-f\circ \gamma =g \; \;  \textrm{où} \; \;  f\in C^{\infty
}(\mathbb{T}^p) \; \; \textrm{est inconnue et} \; \;  g\in
C^{\infty }(\mathbb{T}^p) \; \;  \textrm{est donnée}
\leqno{\hbox{\bf $(\ast )$}}$$

\noindent  Plus précisément, nous montrons que l'image
$\displaystyle Im(\delta )$ de l'opérateur cobord
$$\displaystyle  \delta : C^{\infty }(\mathbb{T}^p)\rightarrow
C^{\infty }(\mathbb{T}^p), \; h\mapsto \delta (h)=h-h\circ \gamma
$$ est un fermé de l'espace de Fréchet $C^{\infty
}(\mathbb{T}^p)$ et que par conséquence l'espace de cohomologie
$$H^1\left(\gamma , C^{\infty }(\mathbb{T}^p)\right):=C^{\infty
}(\mathbb{T}^p)/Im(\delta )$$

\noindent est un espace de Fréchet non trivial. Nous prouvons
également l'existence d'un opérateur linéaire continu $\;
\displaystyle L:Im(\delta )\rightarrow C^{\infty }(\mathbb{T}^p)\;
$ tel que, pour toute fonction $\; \displaystyle g\in Im(\delta ),
\; $  la fonction $\; \displaystyle f=L(g)\; $ est une solution de
l'équation cohomologique {\bf $(\ast )$}. Autrement dit : $\delta
\circ L=id_{Im(\delta )}$.

\vskip 0.2cm

Ce résultat généralise le théorème 5.1 de [DE] prouvé par A.
Dehghan-Nezhad et A. El Kacimi dans le cas où $b=0$ et $A$ est une
matrice diagonalisable dont toutes les valeurs propres sont
réelles et positives.}

\vskip 0.2cm

\vskip 0.5cm

\noindent {\bf 0. Introduction}

\vskip 0.2cm

Si $M$  est une variété différentiable connexe, l'espace vectoriel
$E=C^{\infty }(M)$ des fonctions de classe $C^{\infty }$ sur $ M$
est un espace de Frechet pour la topologie $C^{\infty }$
(topologie de la convergence uniforme de toutes les dérivées sur
les compacts). Une action différentiable d'un groupe dénombrable
$\Gamma $ (supposé de présentation finie) sur $M$ induit sur $E$
une action naturelle donnée par :

$$\forall \gamma \in \Gamma , \; \forall f\in E, \; \gamma .f=f\circ \gamma ^{-1}$$

\vskip 0.2cm

\noindent  On peut donc considérer l'espace vectoriel $H^1(\Gamma
, E)$ désignant le premier espace de cohomologie du groupe $\Gamma
$ à valeurs dans le $\Gamma$-module $E$.

\vskip 0.2cm

 Dans le cas où $\Gamma $  est fini et le
cas où $\Gamma $ agit librement et proprement sur $ M$, (cas où la
projection canonique $\displaystyle M\rightarrow M/\Gamma $ est un
revétement) on montre que $ H^1(\Gamma , E)$ est trivial.

\vskip 0.2cm

Dans le cas où le groupe $\Gamma $  est engendré par un seul
élément $\gamma $, nous avons $\displaystyle \Gamma = \{\gamma ^k
\; | \; k\in \Z\}$  et on peut montrer facilement que l'espace
$H^1(\Gamma , E)$ qu'on notera $H^1(\gamma , E)$ est le conoyau de
l'opérateur cobord :

$$\delta : E\rightarrow E, \; f\mapsto \delta (f)=f-f\circ \gamma $$

Le calcul du conoyau $E/\delta (E)$ de $\delta $  revient à
résoudre l'équation suivante :

$$f-f\circ \gamma =g \hskip 0.5cm \textrm{où} \; \; \cases {
    f\in E \; \; \textrm{est inconnue}\cr
    \textrm{et}\cr
    g\in E \; \; \textrm{est donnée.}}
$$

\noindent dite équation cohomologique associée au système
dynamique discret $\; \gamma :M\rightarrow M$. Les références [A],
[E], [MMY] et [KR] permettent d'avoir une idée sur ce que peut
représenter cette équation dans certains domaines des
mathématiques.

\vskip 0.2cm

{\it Dans ce qui suit, nous allons étudier cette équation dans le
cas d'un automorphisme affine hyperbolique sur un tore [BS].}

\vskip 0.2cm

Soit $A$ est une matrice de $GL(p,\Z )$ qui n'a aucune valeur
propre de module égal à $1$ dans $\C $ et $b$ un élément de
$\R^p$. L'automorphisme affine $\gamma :\R^p\rightarrow \R^p, \;
x\mapsto \gamma (x)=Ax+b$  de $\R^p$ induit sur le tore
$\displaystyle \mathbb{T}^p=\R^p/\Z^p$ un automorphisme affine
hyperbolique que l'on note encore $\gamma $.

L'équation cohomologique associée au système dynamique
hyperbolique $\displaystyle \gamma :\mathbb{T}^p\rightarrow
\mathbb{T}^p$ est la suivante :

$$f-f\circ \gamma =g \hskip 0.5cm \textrm{où} \; \; \cases {
    f\in C^{\infty }(\mathbb{T}^p) \; \; \textrm{est inconnue}\cr
    \textrm{et}\cr
    g\in C^{\infty }(\mathbb{T}^p) \; \; \textrm{est donnée.}} \leqno{\hbox{\bf (0.1)}}
$$

\noindent et l'objet de notre travail est d'établir le théorème
suivant :

\vskip 0.2cm

\noindent {\bf 0.2. Théorème principal.}   {\it Avec les
conditions et les notations ci-dessus, nous avons les résultats
suivants :

\vskip 0.2cm

\noindent {\bf (i)} L'image $\displaystyle Im(\delta )$ de
l'opérateur cobord $\displaystyle \delta : C^{\infty
}(\mathbb{T}^p)\rightarrow C^{\infty }(\mathbb{T}^p), \; h\mapsto
\delta (h)=h-h\circ \gamma $ est un fermé de l'espace de Fréchet
$C^{\infty }(\mathbb{T}^p)$ et par conséquence l'espace de
cohomologie $ H^1\left(\gamma , C^{\infty
}(\mathbb{T}^p)\right):=C^{\infty }(\mathbb{T}^p)/Im(\delta )$ est
un espace de Fréchet non trivial.

\vskip 0.2cm

\noindent {\bf (ii)} Il existe un opérateur linéaire continu $\;
\displaystyle L:Im(\delta )\rightarrow C^{\infty }(\mathbb{T}^p)\;
$ tel que, pour tout élément $g$ de $\; \displaystyle Im(\delta
)$, la fonction $\displaystyle f=L(g)$ est une solution de
l'équation cohomologique ${\bf (0.1)}$. Autrement dit :
$\displaystyle \delta \circ L=id_{Im(\delta )}$.}

\vskip 0.2cm

Avant de commencer la preuve de ce théorème, donnons quelques
exemples d'automorphismes hyperboliques et rappelons quelques
notions qui seront utilisées dans la démonstration.

\vskip 0.4cm

\noindent {\bf 1. Exemples de matrices définissant des
automorphismes hyperboliques sur des tores}

\vskip 0.2cm

\noindent {\bf 1.1. Matrice hyperbolique diagonalisable à valeurs
propres strictement positives :}

\vskip 0.2cm

Un exemple bien connu sur le tore $\; \T^2\; $ est l'automorphisme
hyperbolique (dit du chat d'Arnold) défini par la matrice $\;
\displaystyle A=\left(\begin{array}{cc}
1&1\\
1&2
\end{array}\right)\; $  dont les valeurs propres sont : $\; \lambda _1=\frac{3-\sqrt{5}}{2}\; $
et $\; \lambda _2=\frac{3+\sqrt{5}}{2}$.

\noindent $\; A\; $ est une matrice hyperbolique diagonalisable
dont les valeurs propres sont strictement positives. Elle
appartient donc à la famille $\; \mathcal{F}\; $ des  matrices
étudiées dans l'article [DE].

\vskip 0.2cm

\noindent {\bf 1.2. Matrice hyperbolique diagonalisable à valeurs
propres non toutes strictement positives :}

\vskip 0.2cm

La matrice $\; \displaystyle A=\left(\begin{array}{cc}
1&1\\
1&0
\end{array}\right),\; $  dont les valeurs propres sont : $\; \lambda _1=\frac{1-\sqrt{5}}{2}\; $
et $\; \lambda _2=\frac{1+\sqrt{5}}{2},\; $ définit un
automorphisme hyperbolique sur le tore $\; \T^2$. Cette matrice
est diagonalisable mais ses valeurs propres ne sont pas toutes
strictement positives. Elle n'appartient pas donc à la famille $\;
\mathcal{F}$.

\vskip 0.2cm

\noindent {\bf 1.3. Matrice hyperbolique diagonalisable à valeurs
propres non toutes réelles :}

\vskip 0.2cm

Soit la matrice $\displaystyle \; A=\left(\begin{array}{ccc}
1&1&1\\
1&0&0\\
0&1&0
\end{array}\right)
\; $.  Son polynôme caractéristique $\; P_A(X)\; $ est donné par :

 $$P_A(X)=\hbox{det}(A-XI)=-X^3+X^2+X+1\; \; \textrm{avec} \; \; \hbox{det}(A)=P_A(0)=1$$

\'Etudions les variations de la fonction $\; P_A: x\mapsto P_A(x)\;
$  sur $\R$. Sa dérivée se factorise de la manière suivante
:
$$P'_A(x)=-3x^2+2x+1=(x^2+2x+1)-4x^2=(x+1)^2-(2x)^2=(x+1+2x)(x+1-2x)=(3x+1)(1-x)$$
La fonction $\; P_A\; $ induit donc, par restriction, les trois
bijections suivantes :

$$\; P_A:\left ]-\infty , -\frac{1}{3}\right ]\rightarrow \left [\frac{22}{7} , +\infty \right [ \hskip
0.2cm ; \hskip 0.2cm   P_A:\left [-\frac{1}{3} , 1\right
]\rightarrow \left [\frac{22}{7} , 2\right ]   \hskip 0.2cm ;
\hskip 0.2cm P_A:\left[1 , +\infty \right[\rightarrow \left
]-\infty , 2\right]$$

\vskip 0.2cm

Le polynôme $\; P_A\; $  n'admet donc qu'une seule racine réelle
$\; \mu \; $ située entre les valeurs $\; \frac{3}{2}\; $  et $\;
2\; $ car $\; P_A(\frac{3}{2})=\frac{11}{8}>0\; $  et $\;
P_A(2)=-1<0\; $.

\vskip 0.2cm

Soient $\; \lambda =a+ib\; $  et $\; \overline{\lambda }=a-ib\; $
les deux autres racines de $\; P_A\; $  dans $\; \C\; $  (le
polynôme étant à coefficients réels). On a : $\; b\neq 0 \; $ et $
\; 1=\hbox{det}(A)=\mu \lambda \overline{\lambda }=\mu |\lambda |^2\; $.

Nous avons donc trois valeurs propres  $\; \mu , \; \lambda , \;
\overline{\lambda }\; $   telles que $\; \mu >1\; $  et  $\;
|\lambda |=|\overline{\lambda }|=\frac{1}{\sqrt{\mu }}<1\; $.

\vskip 0.2cm

\noindent {\bf 1.4. Matrice hyperbolique non diagonalisable : }

\vskip 0.2cm

Soit le polynôme $\; Q(X)=\left(P_A(X)\right)^2\; $  où $\;
P_A(X)\; $ est le polynôme de l'exemple 1.3. ci-dessus. $\;
Q(X)=X^6-2X^5-X^4+3X^2+2X+1\in \Z[X]\; $ est à la fois le polynôme
caractéristique et le plynôme minimal de sa matrice compagnon  :

$$C_Q=\left(\begin{array}{cccccc}
0&0&0&0&0&-1\\
1&0&0&0&0&-2\\
0&1&0&0&0&-3\\
0&0&1&0&0&0\\
0&0&0&1&0&1\\
0&0&0&0&1&2
\end{array}\right)$$

\noindent $C_Q\; $ est donc une matrice de $\hbox{GL}(6,\Z)\; $
admettant $\; \mu , \; \lambda \; $ et $\; \overline{\lambda }\; $
comme valeurs propres doubles. Ce qui montre que $\; C_Q\; $ est
une matrice hyperbolique non diagonalisable avec des valeurs
propres qui ne sont pas toutes réelles.

\vskip 0.4cm

\noindent {\bf 2.  L'espace  $\; C^{\infty }(\mathbb{T}^p) \;$}

\vskip 0.2cm

On munit l'espace vectoriel $\; \R^p\; $ de son produit scalaire
usuel et on désigne par $\; C_{per}^{\infty }(\R^p)\; $ l'espace
vectoriel des fonctions $\; \varphi :\R^p\longrightarrow \C \; $
qui sont de classe $\; C^{\infty }\; $ et $\; \Z^p$-péroidiques,
i.e.

$$\forall m\in \Z^p, \; \forall x\in \R^p : \hbox{$\varphi (x+m)=\varphi (x)$}.$$

\noindent Si $\; \pi :\R^p\rightarrow \T^p=\R^p/\Z^p\; $ est la projection
canonique de $\; \R^p\; $ sur le tore $\T^p$, nous avons la
bijection :

$$C^{\infty }(\mathbb{T}^p)\longrightarrow C_{per}^{\infty }(\R^p), \;  h\longmapsto h\circ \pi $$

\noindent permettant d'identifier les éléments de $\; C^{\infty
}(\mathbb{T}^p)\; $ à ceux de $\; C_{per}^{\infty }(\R^p)$. On
utilisera cette identification dans tout ce qui suit. Pour plus de
détail, on pourra consulter la référence [C] par exemple.

\vskip 0.2cm

D'autre part, toute fonction $\; h\in C^{\infty }(\mathbb{T}^p) \;
$ peut être développée en série de Fourier sous la forme :

$$h=\sum _{m\in \Z^p}\widehat{h}(m)\Theta _m$$

\noindent où pour tout $\; m\in \Z^p, \; \widehat{h}(m)\; $  est
le $\; m$-ième coefficient de Fourier de $\; h\; $ donné par :

$$\widehat{h}(m)=\int
_{{\T}^p}h(x)e^{-i2\pi <x,m>} dx=\int _{[0,1]^p}h(t)e^{-i2\pi
<t,m>} dt_1...dt_p$$

\noindent et $\; \Theta _m:\R^p\longrightarrow \C\; $ est la
fonction $\; \Z^p$-périodique définie par :

$$\forall x\in \R^p, \; \hbox{$\Theta _m(x)=e^{i2\pi
<x,m>}$} \hskip 0.5cm   <.,.> \; \; \textrm{étant le produit
scalaire usuel sur}\; \R^p .$$

\noindent On sait que la famille de nombres complexes $\;
\left(\widehat{h}(m)\right)_{m\in \Z^p}\; $ appartient à l'espace
vectoriel $\; \mathcal{S}(\Z^p , \C)\; $  constitué des familles
$\; \displaystyle (a_m)_{m\in \Z^p}$ à décroissance rapide c'est-à-dire telles que :

$$\forall r\in \N, \; \lim _{|m|\rightarrow +\infty }|m|^r|a_m|=0$$

\noindent où pour tout $\; x=(x_1, ..., x_p)\in \R^p$, on $|x|=|x_1|+...+|x_p|$. De plus, l'application :

$$ C^{\infty }(\mathbb{T}^p)\longrightarrow \mathcal{S}(\Z^p , \C), \;
h\longmapsto \left(\widehat{h}(m)\right)_{m\in \Z^p}$$

\noindent est une bijection linéaire qui permet d'identifier $\;
C^{\infty }(\mathbb{T}^p)\; $ à l'espace $\; \mathcal{S}(\Z^p ,
\C)\; $ qu'on peut exprimer sous les deux formes suivantes :

$$\displaystyle \mathcal{S}(\Z^p , \C)=\bigcap _{r\in \N}W^{1,r} \; \; \; \textrm{et}\; \; \; \mathcal{S}(\Z^p , \C)=\bigcap _{r\in \N}W^{2,r}$$

\noindent où pour tout $\; r\in \N, \; W^{1,r} \; $  (resp. $\;
W^{2,r} \; $) désigne l'espace vectoriel des suites de nombres
complexes $\; \displaystyle (a_m)_{m\in \Z^p}\; $ telles que $\;
\displaystyle \sum _{m\in \Z^p} |m|^r|a_m| <\infty \; $  (resp.
$\; \displaystyle \sum _{m\in \Z^p} |m|^{2r}|a_m|^2 <\infty \; $)
muni de la norme

$$\; \displaystyle \|(a_m)\|_{1,r}=|a_0|+\sum _{m\in \Z^p\setminus
\{0\}} |m|^r|a_m| \; \; \;   \left(resp. \; \; \;
\|(a_m)\|_{2,r}=\sqrt{|a_0|+\sum _{m\in \Z^p\setminus \{0\}}
|m|^{2r}|a_m|^2} \right)$$

De plus,  pour tout $\; r\in \N,\; $ les espaces  $\; \; W^{1,r}
\; $ et $\; \; W^{2,r} \; $ sont complets et les injections
suivantes sont des opérateurs compacts :

$$j_{1,r}:W^{1,r+1}\hookrightarrow W^{1,r} \; \; \; \textrm{et}\; \; \; j_{2,r}:W^{2,r+1}\hookrightarrow W^{2,r}$$

\vskip 0.4cm

\noindent {\bf 3. Conditions nécessaires de résolution de
l'équation (0.1)}

\vskip 0.2cm

Pour tout  $k\in \Z$, on note :

 $$\gamma ^k=\cases {
    id \hskip 3.3cm\textrm{si} \; \; k=0 \cr
    \underbrace{\gamma \circ ...\circ \gamma }_{k \; fois} \hskip 2.1cm \textrm{si}\; \;  k>0\cr
    \underbrace{\gamma ^{-1}\circ ...\circ \gamma ^{-1}}_{|k| \; fois} \hskip 1.3cm \textrm{si} \; \; k<0}
$$

\noindent Supposons que l'équation $\; {\bf (0.1)}\; $ admet une
solution $f$. Alors on doit avoir :
$$\int _{{\T}^p} g(x) dx =0$$

\noindent En effet, on a :

$$\int _{{\T}^p}
g(x) dx =\int _{{\T}^p} f(x) dx - \int _{{\T}^p} f\left[\gamma
(x)\right] dx = \int _{{\T}^p} f(x) dx - \int _{{\T}^p} f(t) dt =
0.$$

\vskip 0.2cm

Plus généralement,  nous avons $\; \Phi (g)=0$ pour toute forme
linéaire $\; \Phi :C^{\infty }(\mathbb{T}^p)\rightarrow \C\; $
vérifiant : $\forall h\in C^{\infty }(\mathbb{T}^p), \; \Phi
(h\circ \gamma )=\Phi (h),\; $ .   Autrement dit,

$${\it g \; \; \textrm{annule tout forme linéaire} \; \;  \gamma
-\textrm{invariant sur} \; \; C^{\infty }(\mathbb{T}^p)}
\leqno{\hbox{\bf (3.1)}} $$

\vskip 0.2cm

La prmières condition nécessaire ci-dessus étant un cas
particulier de ${\bf (3.1)}$ puisque la forme linéaire  $\;
h\longmapsto \int _{{\T}^p} h(x) dx\; $ est $\gamma $-invariantes
sur $\; C^{\infty }(\mathbb{T}^p)\; $.

\vskip 0.2cm

La propriété ${\bf (3.1)}$  ci-dessus est donc une condition
nécessaire pour que l'équation cohomologique ${\bf (0.1)}$ admette
au moins une solution pour la donnée $g$.

\vskip 0.2cm

Remarquons que, pour tout $\; k\in \Z,\; $ l'équation $\;
\displaystyle f\circ \gamma ^k - f\circ \gamma ^{k+1} = g\circ
\gamma ^k \; $  est équivalente à l'équation cohomologique ${\bf
(0.1)}$. Ce qui se traduit, au niveau des coefficients de Fourier,
par :

$$\forall k\in \Z,  \; \forall m\in \Z^p, \; \widehat{f\circ \gamma ^k}(m)-\widehat{f\circ \gamma ^{k+1}}(m)=\widehat{g\circ \gamma ^k}(m) \leqno{{\bf (3.2)}}$$

\vskip 0.2cm

Considérons la matrice $B\in \hbox{GL}(p,\Z)$ transposée de la
matrice inverse de $A$ et notons, pour tout $k\in \Z$, $b_k=\gamma
^k(0)$.  On a alors le lemme suivant :

\vskip 0.2cm

\noindent {\bf 3.3. Lemme}  :

\vskip 0.2cm

\noindent {\bf (a)} Pour $k\in \Z$,  $m\in \Z^p$  et  $ h\in
C^{\infty }(\mathbb{T}^p)$, nous avons :
$$b_{-k}=-A^{-k}b_k \; \; \textrm{et} \; \; \widehat{h\circ \gamma ^k}(m)=e^{i2\pi <b_{k},B^km>}\widehat{h}(B^km)$$

\vskip 0.2cm

\noindent {\bf (b)} Pour $ m\in \Z^p\setminus \{0\}$ et $h\in
C^{\infty }(\mathbb{T}^p)$, les trois séries :

$$\displaystyle \sum _{k\geq 0}\widehat{h\circ \gamma ^k}(m), \hskip0.4cm \sum_{k<0}\widehat{h\circ \gamma ^k}(m)=-\sum
_{k\geq 0}\widehat{h\circ \gamma ^{-k}}(m)\hskip0.4cm \hbox{et}
\hskip0.3cm \displaystyle \sum _{k\in \Z}\widehat{h\circ \gamma
^k}(m)$$

\noindent sont absolument convergentes.   De plus,  pour chaque
$m\in \Z^p\setminus \{0\}$, les formes linéaires :

$$\Phi _0:C^{\infty }(\mathbb{T}^p)\rightarrow
\C, \; h\mapsto \Phi _0(h)=\int _{\mathbb{T}^p}h(x)\; dx$$

$$\Phi _m:C^{\infty }(\mathbb{T}^p)\rightarrow
\C, \; h\mapsto \Phi _m(h)=\sum _{k\in \Z}\widehat{h\circ \gamma
^k}(m)$$

$$\Phi ^+_m:C^{\infty }(\mathbb{T}^p)\rightarrow
\C, \; h\mapsto \Phi ^+_m(h)=\sum _{k\geq 0}\widehat{h\circ \gamma
^k}(m)$$

$$\Phi ^-_m:C^{\infty }(\mathbb{T}^p)\rightarrow \C, \; h\mapsto
\Phi ^-_m(h)=-\sum _{k<0}\widehat{h\circ \gamma ^k}(m)$$

\noindent sont continues.

\vskip 0.2cm

\noindent {\bf (c)} Si l'équation cohomologique ${\bf (0.1)}$
admet une solution $f$ pour la donnée $g$ alors :

$$\Phi _0(g)=0, \hskip 0.3cm \Phi _m(g)=0 \hskip 0.3cm \textrm{et} \hskip0.3cm \Phi ^+_m(g)=\widehat{f}(m)=-\Phi ^-_m(g) \hskip0.3cm \textrm{pour tout}\hskip
0.3cm m\in \Z^p\setminus \{0\}$$

\vskip 0.2cm

\noindent {\bf Preuve du Lemme} :

\vskip 0.2cm

\noindent {\bf (a)} Nous avons : \noindent $0=\gamma
^{-k}\left[\gamma ^k(0)\right]=\gamma
^{-k}(b_k)=A^{-k}b_k+b_{-k}$. D'où $b_{-k}=-A^{-k}b_k$.

\vskip 0.2cm

\noindent Nous avons également :
$\displaystyle \widehat{h\circ
\gamma ^k}(m)=\int _{{\T}^p}h\circ \gamma ^k(x)e^{-i2\pi <x,m>}
dx$

\hskip 5.2cm $\displaystyle =\int _{{\T}^p}h\left[\gamma
^k(x)\right]e^{-i2\pi <x,m>} dx$

\hskip 5.2cm $\displaystyle =\int _{{\T}^p}h(u)e^{-i2\pi <\gamma
^{-k}(u),m>} du$

\hskip 5.2cm $\displaystyle =\int _{{\T}^p}h(u)e^{-i2\pi
<A^{-k}(u)+b_{-k},m>} du$

\hskip 5.2cm $\displaystyle =e^{-i2\pi <b_{-k},m>}\int
_{{\T}^p}h(u)e^{-i2\pi <A^{-k}(u),m>} du$

\hskip 5.2cm $\displaystyle =e^{i2\pi <A^{-k}b_{k},m>}\int
_{{\T}^p}h(u)e^{-i2\pi <A^{-k}(u),m>} du$

\hskip 5.2cm $\displaystyle =e^{i2\pi <b_{k},B^km>}\int
_{{\T}^p}h(u)e^{-i2\pi <u,B^km>} du$

\hskip 5.2cm $\displaystyle =e^{i2\pi
<b_{k},B^km>}\widehat{h}(B^km)$

\vskip 0.2cm

\noindent {\bf (b)}   les sommes partielles de chacune des deux
séries $\; \displaystyle \sum _{k\geq 0}|\widehat{h}(B^km)|\; $ et
$\; \displaystyle \sum _{k<0}|\widehat{h}(B^km)|\; $ sont majorées
par la somme $\; \displaystyle \sum _{\alpha \in
\Z^p}|\widehat{h}(\alpha )|\; $. Les trois séries :

$$\displaystyle \sum _{k\geq 0}\widehat{h\circ \gamma ^k}(m)=\sum
_{k\geq 0}e^{i2\pi <b_{k},B^km>}\widehat{h}(B^km)$$

$$\displaystyle \sum _{k<0}\widehat{h\circ \gamma
^k}(m)=\sum_{k<0}e^{i2\pi <b_{k},B^km>}\widehat{h}(B^km)$$

 $$\displaystyle \sum _{k\in \Z}\widehat{h\circ
\gamma ^k}(m)=\sum _{k\in \Z}e^{i2\pi
<b_{k},B^km>}\widehat{h}(B^km)$$

\noindent  sont donc absolument convergentes.  On en déduit que
les formes linéaires $\Phi _m$, $\Phi ^+_m$ et $\Phi ^-_m$ sont
bien définies pour chaque $m\in \Z^p\setminus \{0\}$. De plus,
nous avons :

$$\forall h\in C^{\infty }(\mathbb{T}^p), \; |\Phi _0(h)|=|\widehat{h}(0)|\leq \sum _{\alpha \in
\Z^p}|\widehat{h}(\alpha )|=\|h\|_{1,0}$$

$$\forall h\in C^{\infty }(\mathbb{T}^p), \; |\Phi _m(h)|=|\sum _{k\in \Z}e^{i2\pi
<b_{k},B^km>}\widehat{h}(B^km)|\leq \sum _{k\in
\Z}|\widehat{h}(B^km)|\leq \sum _{\alpha \in
\Z^p}|\widehat{h}(\alpha )|=\|h\|_{1,0}$$

\noindent De même, nous avons : $\displaystyle \forall h\in
C^{\infty }(\mathbb{T}^p), \; |\Phi ^+_m(h)|\leq \|h\|_{1,0}$ et
$\displaystyle |\Phi ^-_m(h)|\leq \|h\|_{1,0}$. Ce qui prouve la
continuité des quatre formes linéaires.

\vskip 0.2cm

\noindent {\bf (c)}  Si $f$ est une solution de ${\bf (0.1)}$ pour
la donnée $g$, alors

$$\Phi _0(g)=\int _{\mathbb{T}^p}g(x) \; dx=\int _{\mathbb{T}^p}f(x)-\int _{\mathbb{T}^p}f\circ \gamma (x) \; dx=0$$

\noindent et si $m\in \Z^p\setminus \{0\}$ et $n\in \N^*$, nous
avons :

$$\sum _{k=0}^n\widehat{g\circ \gamma ^k}(m)=\sum _{k=0}^n\widehat{f\circ \gamma ^k}(m)-\sum _{k=0}^n\widehat{f\circ \gamma ^{k+1}}(m)=
\widehat{f}(m)-\widehat{f\circ \gamma ^{n+1}}(m)$$

\noindent De plus, la série numérique $\displaystyle \sum _{k\geq
0}\widehat{f\circ \gamma ^k}(m)$ est convergente. Donc
$\displaystyle \lim _{n\rightarrow +\infty }\widehat{f\circ \gamma
^{n+1}}(m)=0$ et par suite  $\displaystyle \sum _{k\geq
0}\widehat{g\circ \gamma ^k}(m)=\widehat{f}(m)$ soit
$\displaystyle \Phi ^+_m(g)=\widehat{f}(m)$.

De même,  $\displaystyle \Phi ^-_m(g)=-\sum_{k<0}\widehat{g\circ
\gamma ^k}(m)=\widehat{f}(m)$ et donc aussi $\displaystyle
\displaystyle \sum _{k\in \Z}\widehat{g\circ \gamma ^k}(m)=\Phi
^+_m(g)-\Phi ^-_m(g)=0$.

\hfill $\diamondsuit$

\vskip 0.2cm

Le point (c) du lemme ci-dessus exprime que pour que l'équation
cohomologique ${\bf (0.1)}$ admette une solution $\; f,\; $ pour
la donnée $g$, on doit avoir :

$$\cases{ \displaystyle \forall m\in \Z^p, \; \Phi
_m(g)=0\cr \forall m\in \Z^p\setminus \{0\}, \;
\widehat{f}(m)=\Phi ^+_m(g)=\Phi ^-_m(g).}\leqno{{\bf (3.4)}}$$

\vskip 0.4cm

\noindent {\bf 4. Résultats préliminaires} :

\vskip 0.2cm

\noindent Supposons que la matrice  $\; A\; $ est hyperbolique. La matrice
$\; B,\; $ transposée de l'inverse de $\; A,\; $ est également
hyperbolique puisque ses valeurs propres sont les inverses de
celles de $\; A\; $.

Notons $\; \lambda _1,
...,\lambda _q, \lambda _{q+1}, ..., \lambda _r\; $ les valeurs propres de $\; B\; $
dans $\; \C\; $ avec :
$$\cases{0<|\lambda _j|<1 \hskip0.2cm \textrm{pour}\hskip0.2cm 1\leq j\leq q\cr
|\lambda _j|>1 \hskip0.2cm \textrm{pour}\hskip0.2cm q+1\leq j\leq
r}$$

Pour tout indice $\; j\in \{1, ..., r\},\;  $ on considère
l'espace caractéristique $\; E_j=ker\left(B-\lambda
_jI_p\right)^{k_j}\; $ où $\; k_j\in \N^*\; $  est la multiplicité
de la valeur propre $\; \lambda _j\; $ comme racine du polynôme
minimal $\; P_B(X)\; $. Le lemme des noyaux donne les
décompositions  :

$$\C^p=E_1\oplus ...\oplus E_q\oplus E_{q+1}\oplus ...\oplus E_r \hskip0.5cm \textrm{et}\hskip0.5cm \C^p=E_-\oplus E_+$$

\noindent où $\; E_-=E_1\oplus ...\oplus E_q\; $ (sous-espace
stable) et $\; E_+=E_{q+1}\oplus ...\oplus E_r\; $ (sous-espace
instable). Ces décompositions de l'espace vectoriel $\C^p$
en des sous-espaces invariants par $\; B, \; $ permettent de
définir les projecteurs naturels :

$$\Pi _j:\C^p\longrightarrow \C^p, \; x=x_1+...+x_r\longmapsto x_j
\hskip0.2cm \textrm{pour} \hskip0.2cm 1\leq j\leq r$$

$$\Pi _-:\C^p\longrightarrow \C^p, \; x=x_-+x_+\longmapsto x_-
\hskip0.2cm \textrm{et} \hskip0.2cm \Pi _+:\C^p\longrightarrow
\C^p, \; x=x_-+x_+\longmapsto x_+$$

\noindent ainsi que les automorphismes induits par $\; B\; $ :

$$B_-:E_-\rightarrow E_-, \; x\mapsto B_-(x)=B(x) \; \; \textrm{et}\; \; B_+:E_+\rightarrow E_+, \; x\mapsto B_+(x)=B(x)$$

\vskip 0.2cm

D'autre part, pour tout $\; j\in \{1, ...,r\},\; $ l'opérateur $\;
D_j=\lambda _j\Pi _j \; $ est diagonalisable, l'opérateur $\;
N_j=(B-\lambda _jI_p)\Pi _j\; $ est nilpotent d'indice $\; k_j\; $
et nous avons la décomposition spéctrale :

$$B=N+D\; \; \textrm{où} \; \; \left\{
\begin{array}{lr}
D=D_1+D_2+...+D_r \; \; \textrm{est un opérateur diagonalisable}\\
N=N_1+N_2+...+N_r \; \; \textrm{est un opérateur nilpotent}\\
ND=DN
\end{array}
\right.$$

\noindent Nous avons donc, pour tout entier $\; k,\; $  avec  $\;
\displaystyle  k\geq k_0=\max _{1\leq j\leq r}(k_j),\; $  et pour
$\; x=x_1+...+x_r\in E\; $ :

$\displaystyle B^k(x)=\sum _{j=1}^r B^k(x_j)$

\hskip 1.2cm $\displaystyle =\sum _{j=1}^r\left
[N+D\right]^k(x_j)$

\hskip 1.2cm $\displaystyle =\sum _{j=1}^r\left (\sum _{l=0}^k
C_k^lN^l(x_j)D^{k-l}(x_j)\right)$

\hskip 1.2cm $\displaystyle =\sum _{j=1}^r\left (\sum _{l=0}^k
C_k^lN_j^l(x_j)D_j^{k-l}(x_j)\right)$

\hskip 1.2cm $\displaystyle =\sum _{j=1}^r\left (\sum _{l=0}^k
C_k^l\lambda _j^{k-l}N_j^l(x_j)\right)$

\hskip 1.2cm $\displaystyle =\sum _{j=1}^r \sum _{l=0}^{k_0-1}
C_k^l\lambda _j^{k-l}N_j^l(x_j)$

\vskip 0.2cm

\noindent {\bf 4.1. Remarque}  :  Pour toute valeur
propre $\; \lambda _j\; $ de $\; B,\; \lambda '_j=\frac{1}{\lambda
_j}\; $ est une valeur propre  de $\; B'=B^{-1}\; $  d'ordre $\;
k_j\; $  dont le sous-espace caractéristique $\; E'_j\; $ associé
n'est rien d'autre que $\; E_j$. En effet, nous avons :

$$\left(B-\lambda _jI_p\right)^{k_j}=(-\lambda _j)^{k_j}B^{k_j}\left(B^{-1}-\frac{1}{\lambda
_j}I_p\right)^{k_j}=(-\lambda _j)^{k_j}B^{k_j}\left(B'-\lambda
'_jI_p\right)^{k_j}$$

\noindent Ce qui implique $\; E'_+=E_-, \; E'_-=E_+\; $  et
pour tout  $\; x\in E, \; x'_+=x_-\; $ et $\; x'_-=x_+$.

\vskip 0.2cm

\noindent {\bf 4.2. Lemme}  :  {\it Pour tout $\;
\displaystyle x\in E_-, \; \lim _{k\rightarrow +\infty
}\|B^kx\|=0\; $  et pour tout $\; \displaystyle x\in E_+, \; \lim
_{k\rightarrow +\infty }\|B^{-k}x\|=0$. Par conséquent, nous
avons :  $\; \displaystyle \lim _{k\rightarrow +\infty
}\||B^k_-\||=0\; $  et   $\; \displaystyle \lim _{k\rightarrow
+\infty }\||B^{-k}_+\||=0$,  $\vert \vert \vert \;.\; \vert \vert \vert $  étant la norme induite par le produit
hermitien usuel de $\C^p$ sur chacun des espaces
$\; End(\C ^p), \; End(E_-)\; $ et $\; End(E_+)$.}

\vskip 0.2cm

\noindent {\bf Preuve} :   Soit $\; x=x_1+...+x_q\in
E_-=E_1\oplus ...\oplus E_q$.  Pour $\; k\geq k_0,\; $ nous avons
:

$\displaystyle
\begin{array}{rcl}
\|B^k(x)\|&=&\|\sum _{j=1}^q \sum _{l=0}^{k_0-1} C_k^l\lambda _j^{k-l}N^l(x_j)\|\\[.3cm]
          &\leq &\sum _{j=1}^q \sum _{l=0}^{k_0-1} C_k^l.|\lambda
_j|^{k-l}.\|N^l(x_j)\|
\end{array}
$

\noindent De plus, pour tout $\; j\in \{1, ..., q\},\; $ et pour
tout $\; l\in \{0, ..., k_0-1\},\; $  nous avons $\; \displaystyle
\lim _{k\rightarrow +\infty }C_k^l|\lambda _j|^{k-l}=0\; $  car
$\;|\lambda _j|<1$. Donc   $\; \displaystyle \lim _{k\rightarrow
+\infty }\|B^kx\|=0$. A l'aide de la remarque ci-dessus, on déduit
qu'on a aussi  $\; \displaystyle \lim _{k\rightarrow +\infty
}\|B^{-k}x\|=0\; $ pour $\; x\in E_+$.

\vskip 0.2cm

Pour $\; k\in \N, \; $ nous avons  $\; \displaystyle
\||B^k_-\||=\max \{\|B^kx\| \; | \; x\in E_- \; \textrm{et} \;
\|x\|=1\}=\|B^ku\|\; $  avec $\; u\in E_-$. Donc,  $\;
\displaystyle \lim _{k\rightarrow +\infty }\||B^k_-\||=\lim
_{k\rightarrow +\infty }\|B^ku\|=0$. De même, on a : $\;
\displaystyle \lim _{k\rightarrow +\infty }\||B^{-k}_+\||=0$.

\vskip 0.2cm

\noindent {\bf 4.3. Lemme} : {\it Pour tout $\;
\displaystyle m\in \Z^p\setminus \{0\}, \; $ nous avons : $\;  m _
-\neq 0\; $ et $\; \displaystyle m _+\neq 0$.}

\vskip 0.2cm

\noindent {\bf Preuve} :   Soit $\; m\in \Z^p\setminus
\{0\}$.

Si $\; m_+=0,\; $  alors $\; \displaystyle \lim _{k\rightarrow
+\infty }\|B^k(m)\|^2= \lim _{k\rightarrow +\infty }\|B^k(m
_-)\|^2=0\; $ d'après le lemme (4.2). La suite d'entiers
strictement positifs $\; (\|B^km\|^2)\; $ est donc nulle à partir
d'un certain rang. Ce qui est absurde. D'où $\; m_+\neq 0$. De
même, $\; m_-=m'_+\neq 0\; $ d'après la remarque (4.1). D'où le
lemme.

\vskip 0.2cm

\noindent {\bf 4.4. Proposition}  :  Il existe une
norme $\; \|.\|_*\; $  sur $\; \C ^p=E_-\oplus E_+\; $ telle que :

\vskip 0.2cm

\noindent (a) pour tout $\; (x_- , x_+)\in E_-\times E_+, \;
\|x_-+x_+\|_*=\max (\|x_-\|_* \, , \, \|x_+\|_*)\; $

\vskip 0.2cm

\noindent (b)  $\; \||B_-\||_*<1\; $   et  $\;
\||B^{-1}_+\||_*<1$.

\vskip 0.2cm

On dit que la norme $\; \|.\|_*\; $ est adaptée à l'automorphisme
hyperbolique $B$.

\vskip 0.2cm

\noindent {\bf Preuve}  : D'après le lemme (4.2), $\;
\displaystyle \lim _{k\rightarrow +\infty }\||B^k_-\||=0\; $ et
$\; \displaystyle \lim _{k\rightarrow +\infty }\||B^{-k}_+\||=0$.
Il existe donc un entier $\; n\geq 1\; $ tel que $\;
\||B^n_-\||<1\; $  et $\; \||B^{-n}_+\||<1$. Pour $\; (x_- ,
x_+)\in E_-\times E_+, \; $ posons :

$$\cases { \|x_-\|_*=\sum _{k=0}^{n-1}\|B^kx_-\|\cr
\cr \|x_+\|_*=\sum _{k=0}^{n-1}\|B^{-k}x_+\|\cr \cr
\|x_-+x_+\|_*=\max (\|x_-\|_* , \|x_+\|_*)}$$

Il est clair que $\; \|.\|_*\; $  est une norme sur $\; \C ^p\; $
vérifiant la condition (a). Montrons qu'elle vérifie aussi la
condition (b).

Pour $\; x\in E_-\setminus \{0\}, \; $ nous avons :

$$\frac{\|Bx\|_*}{\|x\|_*}=\frac{\|Bx\|+...\|B^nx\|}{\|x\|+...\|B^{n-1}x\|}=
\frac{\varphi (x)+\frac{\|B^nx\|}{\|x\|}}{\varphi (x)+1}\leq
\frac{\varphi (x)+\||B^n_-\||}{\varphi (x)+1}$$

\noindent où $\; \displaystyle \varphi
(x)=\frac{\|Bx\|}{\|x\|}+...+\frac{\|B^{n-1}x\|}{\|x\|}\leq
\underbrace{\sum _{k=1}^{n-1}\||B^k_-\||}_{\alpha }$.

\vskip 0.2cm

La fonction $\; \displaystyle t\mapsto
\frac{t+\||B^n_-\||}{t+1}=1-\frac{1-\||B^n_-\||}{t+1}\; $ est
strictement croissante sur $\; [0,+\infty [ \; $ et nous avons $\;
0\leq \varphi (x)\leq \alpha \; $. Donc $\; \displaystyle
\frac{\varphi (x)+\||B^n_-\||}{\varphi (x)+1}\leq
 \frac{\alpha +\||B^n_-\||}{\alpha +1}$. D'où :
$\; \displaystyle \||B_-\||_*\leq \frac{\alpha
+\||B^n_-\||}{\alpha +1}<1$.

\vskip 0.2cm

On montre de la même manière que $\;
\||B^{-1}_+\||_*=\||B'_-\||_*<1$.

\vskip 0.2cm

\noindent {\bf 4.5. Remarque}  :  Il est clair que la
norme $\; \|.\|_*\; $ définie ci-dessus  est aussi adaptée à
l'automorphisme hyperbolique $B'=B^{-1}$.

\vskip 0.2cm

\noindent {\bf 4.6. Lemme}  : {\it Soit $\; \displaystyle
m\in \Z^p\setminus \{0\}$. Alors} :

\vskip 0.2cm

\noindent ${\bf (a)}\; $ {\it Si $\; \|m\|_*=\|m_+\|_*, \; $, la
suite $\; (\|B^km\|_*)_{k\geq 0} \; $ est strictement croissante.}

\noindent ${\bf (b)}\; $  {\it Si $\; \|m\|_*=\|m_-\|_*, \; $ alors la
suite $\; (\|B^{-k}m\|_*)_{k\geq 0} \; $ est strictement
croissante.}

\vskip 0.2cm

\noindent {\bf Preuve}  :   Soit $\; m\in \Z^p\setminus
\{0\}$.

\vskip 0.2cm

\noindent ${\bf (a)}\; $ Si $\; \|m\|_*=\|m_+\|_*, \; $ alors

$$\|m\|_*=\|m_+\|_*=\|B^{-1}_+\left(Bm_+\right)\|_*\leq \underbrace{\||B^{-1}_+\||_*}_{<1}.\underbrace{\|Bm_+\|_*}_{\neq 0}<\|Bm_+\|_*\leq \|Bm\|_*$$

$$\|Bm_-\|_*\leq \||B_-\||_*.\|m_-\|_*<\|m_-\|_*\leq \|m\|_*\leq \|Bm_+\|_* \; \; \textrm{et donc} \; \; \|Bm\|_*=\|Bm_+\|_*$$

\noindent Nous avons donc $\; \displaystyle \|m\|_*<\|Bm\|_*\; $
et $\; \displaystyle \|Bm\|_*=\|Bm_+\|_*$. Ce qui permet de
prouver par récurrence que  $\forall k\in \N $, on  a
$\| B^km \|_*<\| B^{k+1}m \|_*$. La suite $\; (\|B^km\|_*)_{k\geq 0} \;
$ est donc strictement croissante.

\vskip 0.2cm

\noindent ${\bf (b)}\; $ Si $\; \|m\|_*=\|m_-\|_*, \; $ alors
d'après les remarques (4.1) et (4.5) nous avons   $\;
\|m\|_*=\|m'_+\|_* \; $  et donc d'après le cas précédent, la
suite $\; (\|B'^km\|_*)_{k\geq 0}=(\|B^{-k}m\|_*)_{k\geq 0} \; $
est strictement croissante.

\vskip 0.4cm

\noindent {\bf 5. Fin de la preuve du théorème}

\vskip 0.2cm

\noindent{\bf 5.1. Preuve du point (i) du théorème} :

\vskip 0.2cm

Pour tout $m\in \Z^p$, la forme linéaire $\displaystyle \Phi _m$
est continue sur $\C^{\infty }(\mathbb{T}^p)$. Ce qui implique que
son noyau $Ker(\Phi _m)$ est un fermé de $C^{\infty
}(\mathbb{T}^p)$ et parsuite l'intersection $\displaystyle
\mathcal{H}:=\cap _{m\in \Z^p}Ker(\Phi _m)$ est un fermé de
$C^{\infty }(\mathbb{T}^p)$.

\noindent D'autre part, si $g\in Im(\delta )$, alors l'équation
${\bf (0.1)}$ admet au moins une solution et donc $g\in
\mathcal{H}$ d'après le point (c) du lemme 3.3. D'où l'inclusion :
$\displaystyle Im(\delta )\subseteq \mathcal{H}$.

\noindent Réciproquement, si $g\in \mathcal{H}$, alors

$$\cases{ \displaystyle \forall m\in \Z^p, \; \Phi _m(g)=0\cr
\forall m\in \Z^p\setminus \{0\}, \;\Phi ^+_m(g)=\Phi ^-_m(g).}$$

\noindent  Montrons que la suite de nombres complexes    $\;
\displaystyle \left(\Phi _m^+(g)\right)_{m\in \Z^p\setminus
\{0\}}\; $ est à décroissance rapide et  définit ainsi une
fonction $\; \displaystyle
 f=\sum _{m\in \Z^p\setminus \{0\}}\Phi _m^+(g)\Theta _m\in C^{\infty
}(\mathbb{T}^p)\; $ solution de  ${\bf (0.1)}$.

\noindent  Pour cela, il suffit de prouver que pour tout $\; r\in
\N,\; $ $\; \displaystyle \lim _{\|m\|_*\rightarrow +\infty
}\|m\|_*^r|\Phi _m^+(g)|=0\; $.

Soit $\; r\in \N\; $ et soit $\;  m\in \Z^p\setminus \{0\}$.

\vskip 0.2cm

\noindent {\bf - Si $\; \|m\|_*=\|m_+\|_* \; $} alors, d'après le
lemme (4.6),  la suite $\; (\|B^km\|_*)_{k\geq 0} \; $ est
strictement croissante et on a donc  :

$\displaystyle \|m\|_*^{r+2}|\Phi
_m^+(g)|=\|m\|_*^{r+2}|\left|\sum _{k\geq 0}\widehat{g\circ \gamma
^k}(m)\right|$

\hskip 2.6cm $\displaystyle \leq \sum _{k\geq
0}\|m\|_*^{r+2}|\widehat{g}(B^km)|$

\hskip 2.6cm $\displaystyle \leq \sum _{k\geq
0}\|B^km\|_*^{r+2}|\widehat{g}(B^km)|$

\hskip 2.6cm $\displaystyle \leq \sum _{\alpha \in \Z^p\setminus
\{0\}}\|\alpha \|_*^{r+2}|\widehat{g}(\alpha )|$

\vskip 0.2cm

\noindent {\bf - Si $\; \|m\|_*=\|m_-\|_* \; $} alors, d'après le
lemme (4.6),  la suite $\; (\|B^{-k}m\|_*)_{k\geq 0} \; $ est
strictement croissante et on a donc  :

$\displaystyle \|m\|_*^{r+2}|\Phi _m^+(g)|=\|m\|_*^{r+2}|\Phi
_m^-(g)|$

\hskip 2.6cm $\displaystyle =\|m\|_*^{r+2}\left|-\sum
_{k<0}\widehat{g\circ \gamma ^k}(m)\right|$

\hskip 2.6cm $\displaystyle =\|m\|_*^{r+2}\left|\sum
_{l>0}\widehat{g\circ \gamma ^{-l}}(m)\right|$

\hskip 2.6cm $\displaystyle \leq \sum
_{l>0}\|m\|_*^{r+2}|\widehat{g}(B^{-l}m)|$

\hskip 2.6cm $\displaystyle \leq \sum
_{l>0}\|B^{-l}m\|_*^{r+2}|\widehat{g}(B^{-l}m)|$

\hskip 2.6cm $\displaystyle \leq \sum _{\alpha \in \Z^p\setminus
\{0\}}\|\alpha \|_*^{r+2}|\widehat{g}(\alpha )|$

\vskip 0.2cm

Les normes   $\displaystyle x=(x_1,...,x_p)\mapsto
|x|=|x_1|+...+|x_p|$ et $\displaystyle x\mapsto \|x\|_*$  étant
équivalentes sur $\C^p$, il existe des nombres réels $\eta
>0$  et $\mu >0$ tels que $\displaystyle \eta \|x\|_*\leq |x|\leq \mu \|x\|_*$  pour tout $x\in
\C^p$. D'où :

$$\displaystyle \forall r\in \N, \; \forall m\in \Z^p\setminus \{0\}, \; \|m\|_*^{r+2}|\Phi _m^+(g)|\leq \frac{\|g\|_{1,r+2}}{\eta
^{r+2}} \; \; \textrm{et}\; \; \|m\|_*^r|\Phi _m^+(g)|\leq
\frac{\|g\|_{1,r+2}}{\eta ^{r+2}\|m\|_*^2}$$

\noindent où $\displaystyle  \|g\|_{1,r+2}=\sum _{\alpha \in
\Z^p\setminus \{0\}}|\alpha |^{r+2}|\widehat{g}(\alpha )|<\infty
$.

\noindent Ce qui implique que $\displaystyle \forall r\in \N, \;
\lim _{\|m\|_*\rightarrow +\infty }\|m\|_*^r|\Phi _m^+(g)|=0$.

Nous avons donc $\; \displaystyle
 f=\sum _{m\in \Z^p\setminus \{0\}}\Phi _m^+(g)\Theta _m\in C^{\infty
}(\mathbb{T}^p)\; $ et d'après le point (a) du lemme 3.3 on a :

\vskip 0.2cm

$\displaystyle \forall m\in \Z^p\setminus\{0\}, \; \widehat{f\circ
\gamma }(m)=e^{i2\pi <b_1,Bm>}\widehat{f}(Bm)$

\vskip 0.2cm

$\displaystyle \hskip 4.2cm=e^{i2\pi <b_1,Bm>}\Phi ^+_{Bm}(g)$

\vskip 0.2cm

$\displaystyle \hskip 4.2cm =\sum _{k\geq 0}e^{i2\pi <b_1,Bm>}
\widehat{g\circ \gamma ^k}(Bm)$

\vskip 0.2cm

$\displaystyle \hskip 4.2cm =\sum _{k\geq 0}\widehat{g\circ \gamma
^{k+1}}(m)$

\noindent D'où : $\displaystyle \forall m\in \Z^p\setminus\{0\},
\; \widehat{\left(f-f\circ \gamma \right)}(m)=\sum _{k\geq
0}\widehat{g\circ \gamma ^k}(m)-\sum _{k\geq 0}\widehat{g\circ
\gamma ^{k+1}}(m)=\widehat{g}(m)$.

\noindent On a aussi $\widehat{\left(f-f\circ \gamma
\right)}(0)=0=\widehat{g}(0)$. Ce qui prouve que $g=\delta (f)\in
Im(\delta )$. On en déduit que $\mathcal{H}\subset Im(\delta )$.

\noindent  $Im(\delta )=\mathcal{H}$  est donc un fermé de
l'espace de Fréchet $C^{\infty }(\T^p)$.  L'espace quotient
$\displaystyle C^{\infty }(\T^p)/Im(\delta )=\mathcal{H}$ qui
n'est rien d'autre que l'espace de cohomologie $H^1\left(\gamma ,
C^{\infty }(\T^p)\right)$ est donc un espace de Fréchet.

\vskip 0.2cm

\noindent{\bf 5.2. Preuve du point (ii) du théorème} :

\vskip 0.2cm

L'opérateur linéaire :

$$L:Im(\delta )\longrightarrow C^{\infty }(\mathbb{T}^p), \; g\longmapsto f=\sum _{m\in \Z^p}\Phi _m^+(g)\Theta _m$$

\noindent est continu. En effet,

$$ \forall r\in \N, \; \forall g\in Im(\delta ), \; \|L(g)\|_{1,r}=\sum
_{m\in \Z^p\setminus \{0\}}|m|^r|\Phi _m^+(g)|\leq
\left[\left(\frac{\mu }{\eta }\right)^{r+2}\sum _{m\in
\Z^p\setminus \{0\}}\frac{1}{|m|^2}\right]\|g\|_{1,r+2}$$

De plus, nous avons : $\displaystyle \forall g\in Im(\delta ), \;
\delta \circ L(g)=g$. \hfill $\diamondsuit $

\bigskip

\centerline{\bf Références}

\bigskip

\noindent [A] {\sc D. V. Anosov} \hskip 0.2cm {\it On an additive
functional homology equation connected with an ergodic rotation of
the circle}. Math. USSR, Izv. 7(1973), 1257-1271.

\vskip 0.2cm

\noindent [C] {\sc Chou Chin-Cheng} \hskip 0.2cm {\it Séries de
Fourier et Théorie des distributions}. Editions Scientifiques
(1983).

\vskip 0.2cm

\noindent [E] {\sc A. El Kacimi Alaoui}  \hskip 0.2cm {\it The
$\overline{\partial }$ operator along the leaves and Guichard's
theorem for a complex simple foliation}. Mathematische Annalen 347
(2010), 885–897.

\vskip 0.2cm

\noindent [DE] {\sc A. Dehghan-Nezhad \& A. El Kacimi Alaoui}
\hskip 0.2cm {\it \'Equations cohomologiques de flots riemanniens
et de diff\'eomorphismes d'Anosov}. Journal of the Mathematical
Society of Japan, Vol. 59 N° 4 (2007), 1105-1134.

\vskip 0.2cm

\noindent [KR] {\sc A. Katok in collaboration with E. A. Robinson,
Jr.} \hskip 0.2cm {\it Cocycles, cohomology and combinatorial
constructions in ergodic theory}. Proceedings of Symposia in Pure
Mathematics Vol. 00 (2001).

\vskip 0.2cm

\noindent [MMY] {\sc S. Marmi, P. Moussa \& J.-C. Yoccoz} \hskip
0.2cm {\it The cohomological equation for roth-type interval
exchange maps}. Journal of the American Mathematical Society, Vol.
18, N° 4 (2005), 823–872.

\vskip 0.2cm

\noindent [BS] {\sc N. Berline \& C. Sabbah} \hskip 0.2cm {\it
Aspects des systèmes dynamiques}. Editions de l'Ecole
Polytechnique, ISBN 978-2-7302-1560-2 (2009).

\vskip 0.4cm

Université Polytechnique Hauts-de-France

ISTV2, LAMAV, FR CNRS 2956

Le Mont Houy

59313 Valenciennes Cedex 9

FRANCE

\vskip 0.2cm

abdellatif.zeggar@uphf.fr

\end{document}